\documentclass[11pt]{article}

\usepackage{amsmath,amssymb,bm}
\usepackage{graphicx}
\usepackage{url}
\usepackage{cite}
\usepackage{theorem}        
\usepackage{enumerate}      

\usepackage[inner=30mm, outer=30mm, textheight=225mm]{geometry}

\newcommand{\co}{\colon\thinspace}

\newcommand{\prooflabel}{Proof}
\newcommand{\qedsymbol}{\rule[-0.5mm]{1.5mm}{3.0mm}}
\newcommand{\qed}{\nolinebreak\hfill\qedsymbol}
\newtheorem{proofthm}{\prooflabel}

\begin{document}

\centerline { \bf A DEFORMATION OF PENNER'S SIMPLICIAL COORDINATE }

\bigskip
\centerline{Tian Yang}
\bigskip
\bigskip

\begin{abstract} We find a one-parameter family of coordinates
$\{\Psi_h\}_{h\in\mathbb{R}}$ which is a deformation of Penner's
simplicial coordinate of the decorated Teichm\"{u}ller space of an ideally
triangulated punctured surface $(S,T)$ of negative Euler
characteristic. If $h\geqslant0$, the decorated Teichm\"{u}ller
space in the $\Psi_h$ coordinate becomes an explicit convex polytope
$P(T)$ independent of $h$, and if $h<0$, the decorated
Teichm\"{u}ller space becomes an explicit bounded convex polytope
$P_h(T)$ so that $P_h(T)\subset P_{h'}(T)$ if $h<h'$. As a
consequence, Bowditch-Epstein and Penner's cell decomposition of the
decorated Teichm\"{u}ller space is reproduced.
\end{abstract}


\bigskip
\centerline{\bf 1. Introduction}
\bigskip

\noindent Decorated Teichm\"{u}ller space of a punctured surface was
introduced by Penner in [P1] as a fiber bundle over the
Teichm\"{u}ller space of complete hyperbolic metrics with cusp ends.
He also gave a cell decomposition of the decorated Teichm\"{u}ller
space invariant under the mapping class group action. To give the
cell decomposition, Penner used the convex hull construction and
introduced the simplicial coordinate $\Psi$ in which the cells can be easily
described. In [BE], Bowditch-Epstein obtained the same cell
decomposition using the Delaunay construction. The corresponding
results for the Teichm\"{u}ller space of a surface with geodesic
boundary have also been obtained. Using Penner's convex hull
construction, Ushijima\,[U] found a mapping class group invariant
cell decomposition, and following the approach of
Bowditch-Epstein\,[BE], Hazel\,[Ha] obtained a natural cell
decomposition of the Teichm\"{u}ller space of a surface with fixed
geodesic boundary lengths. As a counterpart of Penner's simplical
coordinate $\Psi$, Luo\,[L1] introduced a coordinate $\Psi_0$ on the
Teichm\"{u}ller space of an ideally triangulated surface with
geodesic boundary, and Mondello\,[M] pointed out that the $\Psi_0$
coordinate gave a natural cell decomposition of the Teichm\"{u}ller
space.
\\

In [L2], Luo deformed the $\Psi_0$ coordinate to a one-parameter
family of coordinates $\{\Psi_h\}_{h\in\mathbb{R}}$ of the
Teichm\"{u}ller space of a surface with geodesic boundary, and
proved that, for $h\geqslant0$, the image of $\Psi_h$ is an explicit
open convex polytope independent of $h$. For $h<0$, Guo\,[G] proved
that the image of $\Psi_h$ is an explicit bounded open polytope. It
is then a natural question to ask if there is a corresponding deformation of
Penner's simplicial coordinate $\Psi$. The purpose of this paper is to provide
an affirmative answer to this question. We give a one-parameter
family of coordinates $\{\Psi_h\}_{h\in\mathbb{R}}$ of the decorated
Teichm\"{u}ller space of an ideally triangulated punctured surface
so that $\Psi_0$ coincides with Penner's simplicial coordinate $\Psi$\,(Theorem
1.1). We also describe the image of $\Psi_h$\,(Theorem 1.2) and show
that $\Psi_h$ is the unique possible deformation of $\Psi$\,(Theorem 5.1). As an application, Bowditch-Epstein and
Penner's cell decomposition of the decorated Teichm\"{u}ller space
is reproduced using the $\Psi_h$ coordinate\,(Corollary 1.4). The
main results of this paper can be considered as a counterpart of
the work of [G], [L2] and [GL2].
\\

To be precise, let $\overline{T}$ be a triangulation of a closed
surface $\overline{S}$ and let $V$, $E$ and $F$ be the set of
vertices, edges and triangles of $\overline{T}$ respectively. We
call $T=\{\sigma-V\ |\ \sigma\in F\}$ an ideal triangulation of the
punctured surface $S=\overline{S}-V$, and $V$ the set of ideal
vertices (or cusps) of $S$. As a convention in this paper, $S$ is
assumed to have negative Euler characteristic. Let $T_c(S)$ be the
Teichm\"{u}ller space of complete hyperbolic metrics with cusp ends
on $S$. According to Penner\,[P1], a \emph{decorated hyperbolic
metric} $(d,r)\in T_c(S)\times\mathbb{R}_{>0}^V$ on $S$ is the
equivalence class of a hyperbolic metric $d$ in $T_c(S)$ such that
each cusp $v$ is associated with a horodisk $B_v$ centered at $v$ so
that the length of $\partial B_v$ is $r_v$. The space of decorated
hyperbolic metrics $T_c(S)\times\mathbb{R}_{>0}^V$ is the
\emph{decorated Teichm\"{u}ller space}.
\\

Let us recall Penner's simplicial coordinate $\Psi$. Let $(d,r)\in
T_c(S)\times\mathbb{R}_{>0}^V$ be a decorated hyperbolic metric and
let $e$ be an edge of $T$. If $a$ and $a'$ are the generalized
angles (see Section 2) facing $e$, and $b$, $b'$, $c$ and $c'$ are
the generalized angles adjacent to $e$, then Penner's simplicial coordinate
$\Psi\co T_c(S)\times\mathbb{R}_{>0}^V\rightarrow\mathbb{R}^E$ is
defined by

$$\Psi(d,r)(e)=\frac{b+c-a}{2}+\frac{b'+c'-a'}{2}.$$

\begin{figure}[htbp]\centering
\includegraphics[width=5cm]{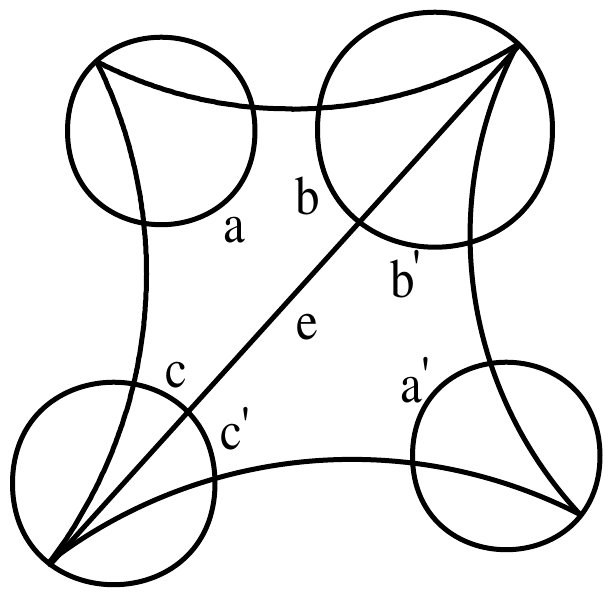}\\
\caption{Penner's simplicial coordinate.}\label{penner}
\end{figure}

An edge path $(t_0,e_1,t_1,\dots,e_n,t_n)$ in a triangulation $T$ is
an alternating sequence of edges $e_i$ with $e_i\neq e_{i+1}$ for $i=1,\dots,n-1$ and triangles $t_i$ so that
adjacent triangles $t_{i-1}$ and $t_i$ share the same edge $e_i$ for
any $i=1,\dots,n$. An \emph{edge loop} is an edge path with
$t_n=t_0$. A \emph{fundamental edge path} is an edge path so that
each edge in the triangulation appears at most once, and a
\emph{fundamental edge loop} is an edge loop so that each edge in
the triangulation appears at most twice. In [P1], Penner proved that
for any vector $z\in\mathbb{R}_{\geqslant0}^E$ such that
$\sum_{i=1}^kz(e_i)>0$ for any fundamental edge loop
$(e_1,t_1,\dots,e_k,t_k)$, there exists a unique decorated complete
hyperbolic metric $(d,r)$ on $S$ so that $\Psi(d,r)=z$. By using a
variational principle on decorated ideal triangles, Guo and
Luo\,[GL1] were able to prove that Penner's simplicial coordinate $\Psi\co
T_c(S)\times\mathbb{R}_{>0}^V\rightarrow\mathbb{R}^E$ is a smooth
embedding with image the convex polytope
$$P(T)=\big\{z\in\mathbb{R}^E\ |\ \sum_{i=1}^kz(e_i)>0\ \text{for any fundamental edge loop}\ (e_1,t_1,\dots,e_k,t_k)\big\}.$$

Let $(S,T)$ be an ideally triangulated punctured surface. To deform Penner's simplicial coordinate, we define for each $h\in\mathbb{R}$ a map $\Psi_h\co
T_c(S)\times\mathbb{R}_{>0}^V\rightarrow\mathbb{R}^E$ by

$$\Psi_h(d,r)(e)=\int_0^{\frac{b+c-a}{2}}e^{ht^2}dt+\int_0^{\frac{b'+c'-a'}{2}}e^{ht^2}dt,$$

\noindent where $a$ and $a'$ are the generalized angles facing $e$,
and $b$, $b'$, $c$ and $c'$ are the generalized angles adjacent to
$e$ as in Figure \ref{penner}. The main theorems of this paper are the following

\medskip
\noindent {\bf Theorem 1.1} \it Suppose that $(S,T)$ is an ideally
triangulated punctured surface. Then for all $h\in\mathbb{R}$, the
map $\Psi_h\co T_c(S)\times\mathbb{R}_{>0}^V\rightarrow\mathbb{R}^E$
is a smooth embedding.\rm
\medskip

\noindent {\bf Theorem 1.2} \it For $h\in\mathbb{R}$ and an ideally
triangulated punctured surface $(S,T)$, let $P_h(T)$ be the set of
points $z\in\mathbb{R}^E$ such that
\begin{enumerate}[(a)]
\item $z(e)<2\int_0^{+\infty}e^{ht^2}dt$ for each edge $e\in E$,
\item $\sum_{i=1}^nz(e_i)>-2\int_0^{+\infty}e^{ht^2}dt$ for each fundamental edge path $(t_0,e_1,t_1,\dots,e_n,t_n)$,
\item $\sum_{i=1}^nz(e_i)>0$ for each fundamental edge loop $(e_1,t_1,\dots,e_n,t_n)$.
\end{enumerate}
Then we have $\Psi_h(T_c(S)\times\mathbb{R}_{>0}^V)=P_h(T)$.
Furthermore, if $h\geqslant0$, then conditions (a) and (b) become
trivial, and the image of $\Psi_h$ is the open convex polytope
$P(T)$, hence independent of $h$; and if $h<0$, then the image
$P_h(T)$ is a bounded open convex polytope so that $P_h(T)\subset
P_{h'}(T)$ if $h<h'$.\rm
\medskip

Clearly $\Psi_0$ coincides with Penner's simplicial coordinate $\Psi$ and
$\Psi_h$ is a deformation of $\Psi$. Theorem 1.1 is
proved in Section 2 using the strategy of Guo-Luo\,[GL1]. We set up
a variational principle from the derivative cosine law of decorated
ideal triangles whose energy function $V_h$ is strictly concave. For
$i=1,\dots,|E|$, each variable of $V_h$ is a smooth monotonic function
of the edge length $l_i$ in the decorated hyperbolic metric $(d,r)$,
and $\Psi_h$ is the gradient of $V_h$, hence is a smooth embedding.
We study various degenerations of decorated ideal triangles in
Section 3 with which we will prove Theorem 1.2 in Section 4. We will
also prove that $\{\Psi_h\}_{h\in\mathbb{R}}$ is the unique possible
deformation of Penner's simplicial coordinate by using a variational principle
(Theorem 5.1).
\\

The Delaunay cell decomposition of a decorated hyperbolic surface
will be reviewed in Section 6 and we will prove the following

\medskip
\noindent {\bf Theorem 1.3} \it Suppose $(S,T)$ is an ideally
triangulated punctured surface, and $(d,r)\in
T_c(S)\times\mathbb{R}_{>0}^V$ is a decorated hyperbolic metric so
that the horodisks associated to the ideal vertices do not
intersect. Then for all $h\in\mathbb{R}$, the corresponding Delaunay decomposition
$\Sigma_{d,r}$ coincides with the ideal triangulation $T$ if and
only if $\Psi_h(d,r)(e)>0$ for each $e\in E$.\rm
\medskip

{Bowditch-Epstein\,[BE] and Penner\,[P1] showed that there is a
natural cell decomposition of the decorated Teichm\"{u}ller space
$T_c(S)\times\mathbb{R}_{>0}^V$ invariant under the mapping class
group action. One interesting consequence of Theorems 1.1, 1.2 and
1.3 is the following. Let $A(S)-A_{\infty}(S)$ be the fillable arc
complex as in [H], and let $|A(S)-A_{\infty}(S)|$ be its underlying
space. Penner\,[P1] provided a mapping class group equivariant
homeomorphism

$$\Pi\co T_c(S)\times\mathbb{R}_{>0}^V\rightarrow
|A(S)-A_{\infty}(S)|\times\mathbb{R}_{>0}$$ \noindent so that the
restriction of $\Pi$ to each simplex of maximum dimension is given
by the simplicial coordinate $\Psi$. Using Penner's method, we have the
following

\medskip
\noindent {\bf Corollary 1.4} \it Suppose that $S$ is a punctured
surface of negative Euler characteristic.

\begin{enumerate}[(a)]
\item For all $h>0$, there is a homeomorphism

$$\Pi_h\co T_c(S)\times\mathbb{R}_{>0}^V\rightarrow
|A(S)-A_{\infty}(S)|\times\mathbb{R}_{>0}$$

equivariant under the mapping class group action so that the
restriction of $\Pi_h$ to each simplex of maximum dimension is given
by the $\Psi_h$ coordinate.

\item The cell structures for various $h>0$ are the same as Penner's.\
\end{enumerate}\rm
\medskip

\bigskip
\noindent{\bf Acknowledgment} The author would like to thank Feng
Luo for instructive discussions on this subject and suggestions for
improving this paper, Shiu-chun Wong for making several crucial
suggestions on polishing the writing of this paper, Ren Guo, Julien Roger and Jian Song for useful suggestions and Tianling
Jin for helpful discussions. The author is also very grateful to the referee for the carefully reading and making many valuable suggestions on both the mathematics and the writing of this paper. Part of the work is supported by an NSF
research fellowship.


\bigskip
\centerline{\bf 2. A variational principle on decorated ideal
triangles}\label{2}
\bigskip

\noindent Let $(S,T)$ be an ideally triangulated punctured surface
with a set of ideal vertices $V$ and a set of edges $E$. We assume
that $S$ has negative Euler characteristic. The proof of Theorem 1.1
goes as follows. By Penner [P1], there is a smooth parametrization
of the decorated Teichm\"{u}ller space
$T_c(S)\times\mathbb{R}_{>0}^V$ by $\mathbb{R}^E$ using the edge
lengths. From the cosine law of decorated ideal triangles\,[P1], we
construct for each $h\in\mathbb{R}$ a smooth strictly convex
function $V_h$ on a convex subset of $\mathbb{R}^E$ so that its
gradient is $\Psi_h$. By a variational principle, for each
$h\in\mathbb{R}$, the map $\Psi\co
T_c(S)\times\mathbb{R}_{>0}^V\rightarrow R^E$ is a smooth embedding.
This variational principle, whose proof is elementary, is: \emph{If
$X$ is an open convex set in $\mathbb{R}^n$ and $f\co
X\rightarrow\mathbb{R}$ is smooth strictly concave, then the
gradient $\nabla f\co X\rightarrow\mathbb{R}^n$ is injective.
Furthermore, if the Hessian of $f$ is negative definite for all
$x\in X$, then $\nabla f$ is a smooth embedding.}
\medskip

A \emph{decorated ideal triangle} $\Delta$ in the hyperbolic plane
$\mathbb{H}^2$ is an ideal triangle such that each ideal vertex $v$
is associated with a horodisk $B_v$ centered at $v$. If $e_1$ and
$e_2$ are two edges adjacent to an ideal vertex $v$ of $\Delta$,
then the \emph{generalized angle} of $\Delta$ at $v$ is defined to
be the length of the intersection of $\partial B_v$ and the cusp
region enclosed by $e_1$ and $e_2$. (In\,[P1], Penner called the generalized angles the $h$-lengths of a decorated ideal triangle, and in\,[GL1], Guo and Luo defined
the generalized angle to be twice of the generalized angle defined
here.) If $e$ is an edge of $\Delta$ with ideal vertices $u$ and
$v$, then the \emph{generalized edge length} (or \emph{edge length}
for simplicity) of $e$ in $\Delta$ is the signed hyperbolic distance
between the intersection of $e$ and $\partial B_u$ and the
intersection of $e$ and $\partial B_v$ (Figure \ref{figure1}\,(a)).
Note that if $B_u\cap B_v\neq\emptyset$, then the generalized edge
length of $e$ is either zero or negative (Figure
\ref{figure1}\,(b)). In a decorated hyperbolic metric $(d,r)\in
T_c(S)\times\mathbb{R}_{>0}^V$, each triangle $\sigma$ in $T$ is
isometric to an ideal triangle and the decoration
$r\in\mathbb{R}_{>0}^V$ induces a decoration on $\sigma$. If $e\in
E$ is an edge and $\sigma$ is an ideal triangle adjacent to $e$,
then the \emph{generalized edge length} $l_{d,r}(e)$ of $e$ is
defined to be the generalized edge length of $e$ in $\sigma$. It is
clear that $l_{d,r}(e)$ does not depend on the choice of $\sigma$.

\begin{figure}[htbp]\centering
\includegraphics[width=10cm]{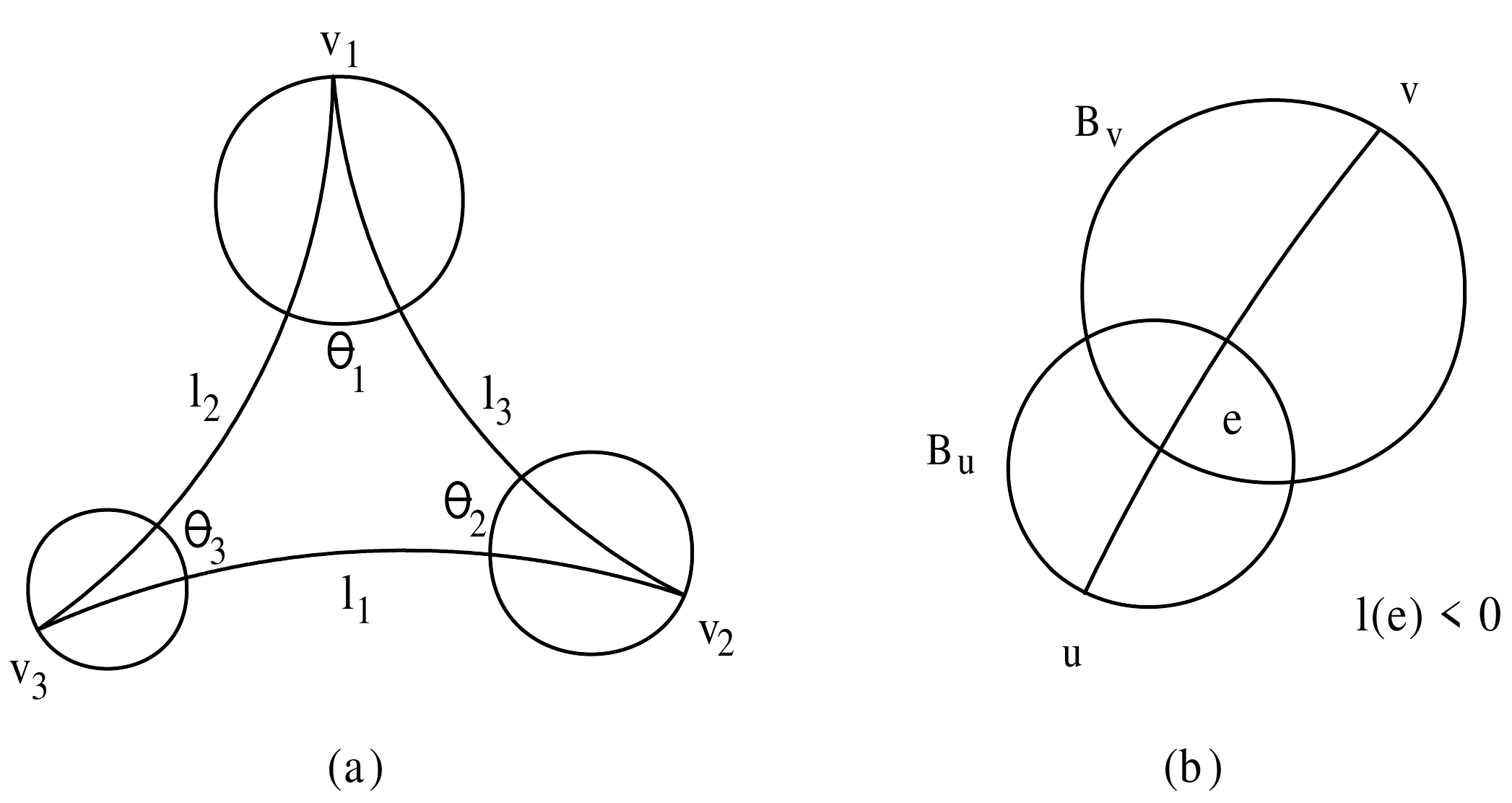}\\
\caption{Generalized angles and edge lengths.}\label{figure1}
\end{figure}

Penner\,[P1] defined the length parametrization
\begin{equation*}
\begin{split}
L\co T_c(S)\times\mathbb{R}_{>0}^V\rightarrow&\mathbb{R}^E\\
(d,r)\mapsto &l_{d,r}
\end{split}
\end{equation*}
and showed that $L$ is a diffeomorphism. (The exponential of half of the generalized edge length, which is called the $\lambda$-length in\,[P1], is sometimes called Penner's coordinate in the literature.) Penner also proved the
following cosine law of decorated ideal triangles. Suppose that
$\Delta$ is a decorated ideal triangle with edge lengths $l_1$,
$l_2$ and $l_3$ and opposite generalized angles $\theta_1$,
$\theta_2$ and $\theta_3$. For $i,j,k=1,2,3$,
\begin{equation}\label{cosinelaw}
\theta_i=e^{\frac{l_i-l_j-l_k}{2}}\ \ \text{and}\ \ e^{l_i}=\frac{1}{\theta_j\theta_k}.
\end{equation}
As a consequence, there is the sine law of decorated triangles:
\begin{equation}\label{sinelaw}
\frac{\theta_1}{e^{l_1}}=\frac{\theta_2}{e^{l_2}}=\frac{\theta_3}{e^{l_3}}.
\end{equation}
\\

For $i,j,k=1,2,3$ and $x_i=\frac{\theta_j+\theta_k-\theta_i}{2}$,
let $\mu(x_i)=\int_0^{x_i}e^{ht^2}dt$ and
$u_i=\int_0^{l_i}e^{-he^{-t}}dt$. Denote by $U\subset\mathbb{R}^3$
the set of all possible values of $u=(u_1,u_2,u_3)$.
\medskip

\noindent {\bf Lemma 2.2} \it For each $h\in\mathbb{R}$, the
differential $1$-form $\omega_h=\sum_{i=1}^3\mu(x_i)du_i$ is closed
in $U$ and the function $F_h$ defined by the integral
$F_h(u)=\int_0^u\omega_h$ is strictly concave in $U$. Furthermore,
\begin{equation*}
\frac{\partial F_h}{\partial u_i}=\int_0^{x_i}e^{ht^2}dt.
\end{equation*}\rm
\medskip

\noindent \textit{Proof}: Consider the matrix $H=[\frac{\partial
\mu(x_i)}{\partial u_j}]_{3\times 3}$. The closedness of $\omega_h$
is equivalent to that $H$ is symmetric, and the strict concavity of
$F_h$ will follow from the negative definiteness of $H$. It follows
from the partial derivatives of (\ref{cosinelaw}) that
$\frac{\partial x_i}{\partial l_i}=-\frac{x_i+x_j+x_k}{2}$ and
$\frac{\partial x_i}{\partial l_j}=\frac{x_k}{2}$. We have

\begin{equation*}
\frac{\partial \mu(x_i)}{\partial u_i}=\frac{e^{hx_i^2}}{e^{-he^{-l_i}}}\frac{\partial x_i}{\partial l_i}=-\frac{x_i+x_j+x_k}{2}e^{h(\frac{\theta_i^2+\theta_j^2+\theta_k^2}{4}+\frac{3\theta_j\theta_k-\theta_i\theta_k-\theta_i\theta_j}{2})},
\end{equation*}
and for $i\neq j$, we have

\begin{equation*}
\frac{\partial \mu(x_i)}{\partial u_j}=\frac{e^{hx_i^2}}{e^{-he^{-l_j}}}\frac{\partial x_i}{\partial l_j}=\frac{x_k}{2}e^{h(\frac{\theta_i^2+\theta_j^2+\theta_k^2}{4}+\frac{\theta_j\theta_k+\theta_i\theta_k-\theta_i\theta_j}{2})},
\end{equation*}
from which we see that $H$ is symmetric. Let
$c=\frac{1}{2}e^{{h(\frac{\theta_i^2+\theta_j^2+\theta_k^2}{4}-\frac{\theta_j\theta_k+\theta_i\theta_k+\theta_i\theta_j}{2})}}>0$
and let $D$ be the diagonal matrix whose $(i,i)$-th entry is
$e^{h\theta_j\theta_k}$. The matrix $H$ can be written as $cDMD$,
where
\begin{equation*}
M=
\begin{bmatrix}
-(x_1+x_2+x_3) & x_3 & x_2\\
x_3 & -(x_1+x_2+x_3) & x_1\\
x_2 & x_1 & -(x_1+x_2+x_3)
\end{bmatrix}.
\end{equation*}
The negative definiteness of $H$ is equivalent to that of $M$, i.e., the positive definiteness of $-M$. This follows from the direct calculation that each leading principal minor is positive using Sylvester's criterion. \hfill{q.e.d}
\\

\noindent \textit{Proof of Theorem 1.1}: For a decorated hyperbolic
metric $(d,r)\in T_c(S)\times\mathbb{R}_{>0}^V$, let
$l_{d,r}\in\mathbb{R}^E$ be its length parameter. The integral $u(e)=\int_0^{l_{d,r}(e)}e^{-he^{-t}}dt$ is a smooth monotonic
function of $l_{d,r}(e)$, and the possible values of $u$ form an
open convex cube $U$ in $\mathbb{R}^E$. With $u_i=u(e_i)$, the
energy function $V_h\co U\rightarrow \mathbb{R}$ is defined by

$$V_h(u)=\sum_{\{e_i,e_j,e_k\}} F_h(u_i,u_j,u_k),$$

\noindent in which the summation is taken over all of the decorated
ideal triangles. By Lemma 2.2, $V_h$ is smooth and strictly concave
in $U$ and

$$\frac{\partial V_h}{\partial u_i}=\Psi_h(e_i),$$

\noindent i.e., $\nabla V_h=\Psi_h$. By the variational principle,
the map $\Psi_h=\nabla V_h\co U\rightarrow \mathbb{R}^E$ is a smooth
embedding. \hfill{q.e.d}


\bigskip
\centerline {\bf 3. Degenerations of decorated ideal triangles}
\bigskip

\noindent To describe the image of $\Psi_h$, we study degenerations
of decorated ideal triangles. Suppose $\Delta$ is a decorated ideal
triangle with edge lengths $l_1$, $l_2$ and $l_3$ and opposite
generalized angles $\theta_1$, $\theta_2$ and $\theta_3$.

\medskip
\noindent {\bf Lemma 3.1}\it \begin{enumerate}[(I)]
\item If $\{(l_1,l_2,l_3)\}$ converges to $(-\infty,c_2,c_3)$ with
$c_2, c_3\in(-\infty,+\infty]$, then $\{\theta_1\}$ converges to
$0$, and we can take a subsequence so that at least one of
$\{\theta_2\}$ and $\{\theta_3\}$ converges to $+\infty$.
\item If $\{(l_1,l_2,l_3)\}$ converges to $(-\infty,-\infty,c_3)$ with
$c_3\in(-\infty,+\infty]$, then $\{\theta_3\}$ converges to
$+\infty$, and we can take a subsequence so that at least one of
$\{\theta_1\}$ and $\{\theta_2\}$ converges to a finite number.
\item If $\{(l_1,l_2,l_3)\}$ converges to $(-\infty,-\infty,-\infty)$,
then we can take a subsequence such that at least two of
$\{\theta_1\}$, $\{\theta_2\}$ and $\{\theta_3\}$ converge to
$+\infty$.
\end{enumerate}\rm
\medskip

\noindent \textit{Proof}: For (I), if $\{(l_1,l_2,l_3)\}$ converges
to $(-\infty,c_2,c_3)$, then $\{\frac{l_1-l_2-l_3}{2}\}$ converges
to $-\infty$. By cosine law (\ref{cosinelaw}),
$\{\theta_1\}=\{e^{\frac{l_1-l_2-l_3}{2}}\}$ converges to $0$. Let
$a_2=\frac{l_2-l_1-l_3}{2}$ and $a_3=\frac{l_3-l_1-l_2}{2}$, so
$\{a_2+a_3\}=\{-l_1\}$ converges to $+\infty$. Thus, by taking a
subsequence if necessary, at least one of $\{a_2\}$ and $\{a_3\}$,
say $\{a_2\}$, converges to $+\infty$, and
$\{\theta_2\}=\{e^{a_2}\}$ converges to $+\infty$. For (II), if
$\{(l_1,l_2,l_3)\}$ converges to $(-\infty,-\infty,c_3)$, then
$\{\frac{l_3-l_1-l_2}{2}\}$ converges to $+\infty$, and
$\{\theta_3\}=\{e^{\frac{l_3-l_1-l_2}{2}}\}$ converges to $+\infty$.
Letting $a_1=\frac{l_1-l_2-l_3}{2}$ and $a_2=\frac{l_2-l_1-l_3}{2}$,
we have $\{a_1+a_2\}=\{-l_3\}$ converges to $-c_3$. Thus, either both
$\{a_1\}$ and $\{a_2\}$ converge to a finite number, or by taking a
subsequence if necessary, at least one of $\{a_1\}$ and $\{a_2\}$,
say $\{a_1\}$, converges to $-\infty$. In the former case, both
$\{\theta_1\}=\{e^{a_1}\}$ and $\{\theta_2\}=\{e^{a_2}\}$ converge
to a finite number, and in the latter case,
$\{\theta_1\}=\{e^{a_1}\}$ converges to $0$. For (III), we have by
cosine law (\ref{cosinelaw}) that
$\{\theta_1\theta_2\}=\{e^{-l_3}\}$ converges to $+\infty$. Thus, by
taking a subsequence if necessary, at least one of $\{\theta_1\}$
and $\{\theta_2\}$, say $\{\theta_1\}$, converges to $+\infty$.
Since $\{\theta_2\theta_3\}=\{e^{-l_1}\}$ converges to $+\infty$ as
well, by taking a subsequence, at least one of $\{\theta_2\}$ and
$\{\theta_3\}$ converges to $+\infty$. \hfill{q.e.d}
\\

We call a converging sequence of decorated ideal triangles in (I),
(II) and (III) of Lemma 3.1 a \emph{degenerated decorated ideal
triangle of type I, II} and \emph{III} respectively. If $a$ is the
generalized angle facing an edge $e$ in a decorated triangle
$\Delta$, and $b$ and $c$ are the generalized angles adjacent to
$e$, then we call $x(e)=\frac{b+c-a}{2}$ the \emph{$x$-invariant} of
$e$ in $\Delta$.

\medskip
\noindent {\bf Corollary 3.2} \it If $\Delta$ is a degenerated
decorated ideal triangle of type I, II or III, then by taking a
subsequence if necessary, there is an edge $e$ of $\Delta$ such that
$\{l(e)\}$ converges to $-\infty$ and $\{x(e)\}$ converges to
$+\infty$.\rm
\medskip

\noindent \textit{Proof}: If $\Delta$ is of type I and $\{l_1\}$
converges to $-\infty$, then by Lemma 3.1 (I),
$\{x_1\}=\{\frac{\theta_2+\theta_3-\theta_1}{2}\}$ converges to
$+\infty$. If $\Delta$ is of type II and $\{(l_1,l_2,l_3)\}$
converges to $(-\infty,-\infty,c_3)$, then by Lemma 3.1 and taking a
subsequence if necessary, at least one of $\{\theta_1\}$ and
$\{\theta_2\}$, say $\{\theta_1\}$, converges to a finite number,
and $\{\theta_3\}$ converges to $+\infty$. Thus, $\{l_1\}$ converges
to $-\infty$ and $\{x_1\}=\{\frac{\theta_2+\theta_3-\theta_1}{2}\}$
converges to $+\infty$. If $\Delta$ is of type III, then there are
at least two of $\{\theta_1\}$, $\{\theta_2\}$ and $\{\theta_3\}$
that converge to $+\infty$. Suppose $\{\theta_3\}$ is one of the two
that converge to $+\infty$. Since $\{x_1+x_2\}=\{\theta_3\}$
converges to $+\infty$, by taking a subsequence if necessary, at
least one of $\{x_1\}$ and $\{x_2\}$, say $\{x_1\}$, converges to
$+\infty$. Thus, $\{l_1\}$ converges to $-\infty$ and $\{x_1\}$
converges to $+\infty$. \hfill{q.e.d}
\\

We call an edge $e$ as in Corollary 3.2 where $l(e)\to-\infty$ and $x(e)\to+\infty$ a \emph{bad edge} of $\Delta$,
and otherwise, $e$ is a \emph{good edge}. Note that there may be more than one bad
edge in a degenerated ideal triangle.

\medskip
\noindent {\bf Lemma 3.3} \it Let $\{\Delta^{(m)}\}$ be a sequence
of decorated ideal triangles that converges to a degenerated
decorated ideal triangle $\Delta$ of type I, II or III. Then we can
take a subsequence so that for $m$ sufficiently large, the length of
each bad edge of $\Delta^{(m)}$ is strictly less than the length of
each good edge.\rm
\medskip

\noindent \textit{Proof}: If $\Delta$ is of type I, then by Lemma
3.1, the length of the only bad edge converges to $-\infty$ and the
length of other two edges converge to a finite number. For $m$
sufficiently large, the length of the bad edge is less than the
lengths of the good edges.
\\
\begin{figure}[htbp]\centering
\includegraphics[width=10cm]{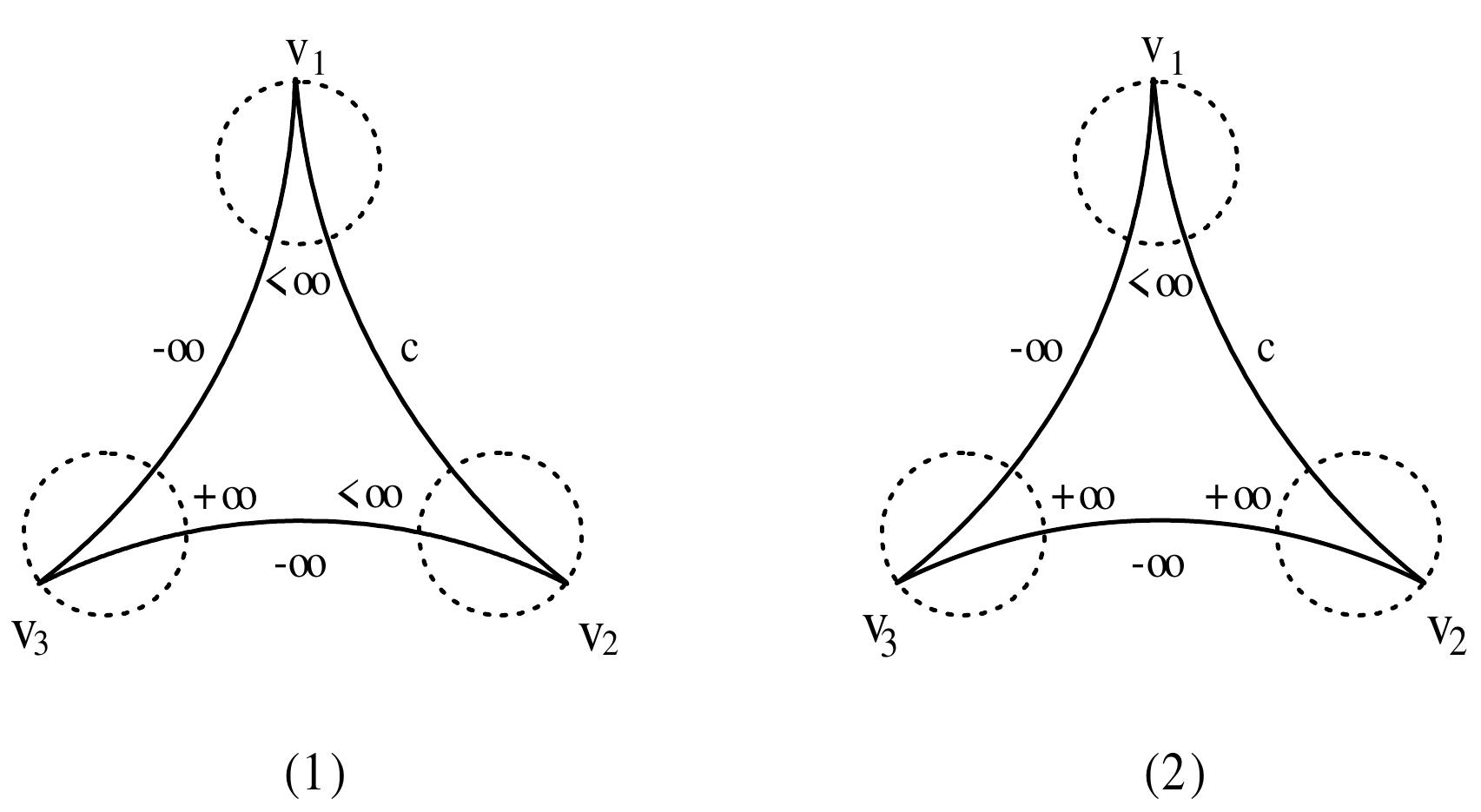}\\
\caption{Type II.}\label{type2}
\end{figure}

If $\Delta$ is of type II, we may assume that
$\{(l_1^{(m)},l_2^{(m)},l_3^{(m)})\}$ converges to
$(-\infty,-\infty,c)$ with $c\in(-\infty,+\infty]$. By Lemma 3.1,
there are two cases to be considered (Figure \ref{type2}).
\\

\noindent \textit{Case 1.} Suppose that $\theta_3^{(m)}$ converges
to $+\infty$ and both $\theta_1^{(m)}$ and $\theta_2^{(m)}$ converge
to a finite number. In this case, both $l_1$ and $l_2$ are bad and
converge to $-\infty$. The only good edge length $l_3$ converges to
$c\in(-\infty,+\infty]$. Hence for $m$ sufficiently large,
$l_1^{(m)}<l_3^{(m)}$ and $l_2^{(m)}<l_3^{(m)}$.
\\

\noindent \textit{Case 2.} Suppose that $\theta_3^{(m)}$ converges
to $+\infty$, and one of $\theta_1^{(m)}$ and $\theta_2^{(m)}$, say
$\theta_2^{(m)}$, converges to $+\infty$ and $\theta_1^{(m)}$
converges to a finite number. In this case $l_1$ is bad. If $l_2$ is
also bad, then both $l_1$ and $l_2$ converge to $-\infty$, and $l_3$
converges to $c\in(-\infty,+\infty]$. Hence for $m$ sufficiently
large, $l_1^{(m)}<l_3^{(m)}$ and $l_2^{(m)}<l_3^{(m)}$. If $l_2$ is
good, then $\theta_1^{(m)}<\theta_2^{(m)}$ for $m$ sufficiently
large, since $\theta_1^{(m)}$ converges to a finite number and
$\theta_2^{(m)}$ converges to $+\infty$. By sine law
(\ref{sinelaw}), $l_1^{(m)}<l_2^{(m)}$.
\\

\begin{figure}[htbp]\centering
\includegraphics[width=10cm]{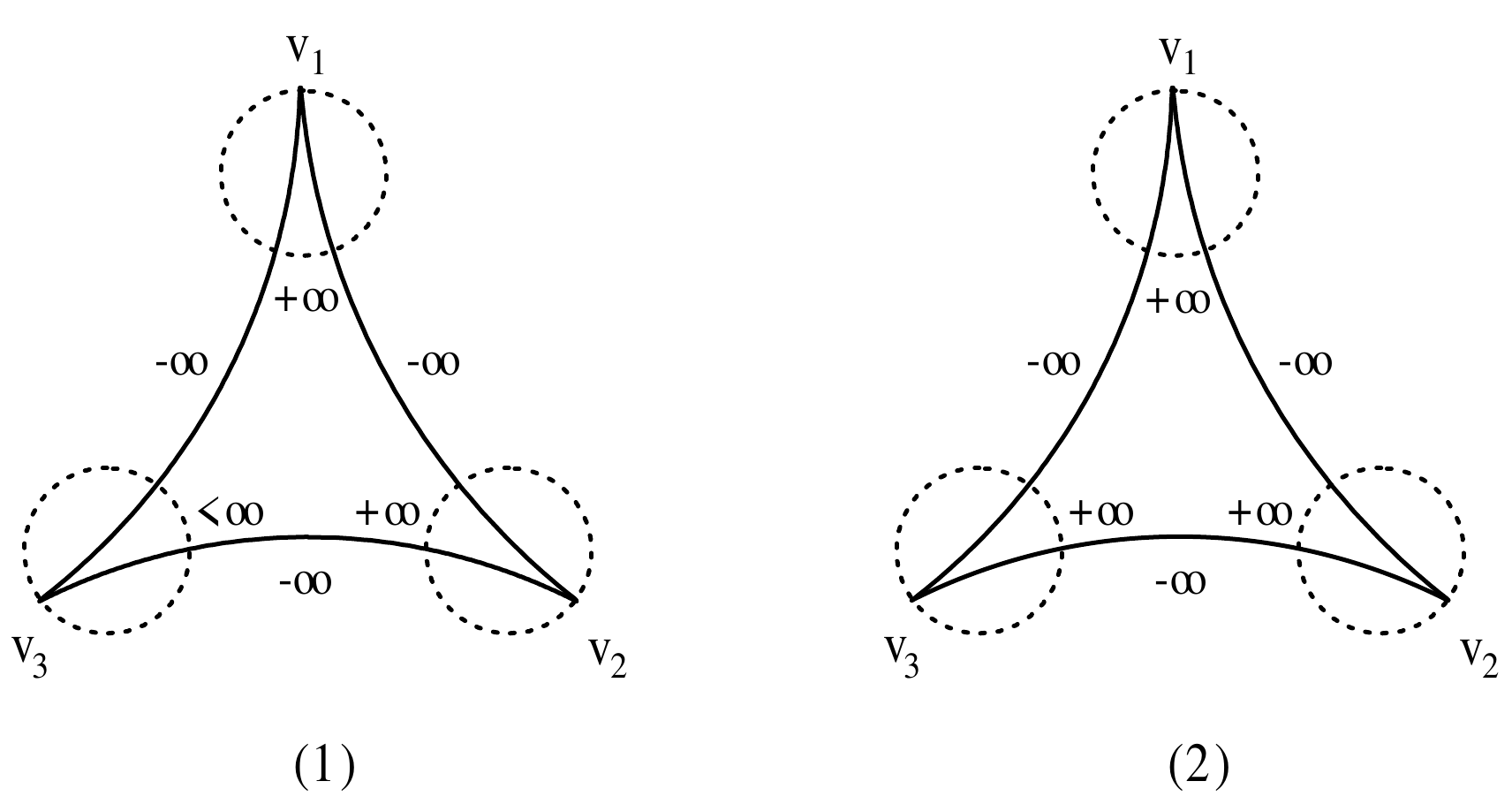}\\
\caption{Type III.}\label{type3}
\end{figure}

If $\Delta$ is of type III, then by Lemma 3.1, we also consider two
cases (Figure \ref{type3}).
\\

\noindent \textit{Case 1.} Two of $\theta_1^{(m)}$, $\theta_2^{(m)}$
and $\theta_3^{(m)}$, say $\theta_1^{(m)}$ and $\theta_2^{(m)}$
converge to $+\infty$, and $\theta_3^{(m)}$ converges to a finite
number. In this case, $l_3$ is bad. Since
$\theta_3^{(m)}<\theta_1^{(m)}$ and $\theta_3^{(m)}<\theta_2^{(m)}$
for $m$ sufficiently large, by sine law (\ref{sinelaw}),
$l_3^{(m)}<l_1^{(m)}$ and $l_3^{(m)}<l_2^{(m)}$. If one of $l_1$ and
$l_2$, say $l_2$, is also bad, then
$x_2^{(m)}=\frac{\theta_1^{(m)}+\theta_3^{(m)}-\theta_2^{(m)}}{2}$
converges to $+\infty$. Since $\theta_3^{(m)}$ converges to a finite
number, $\theta_2^{(m)}<\theta_1^{(m)}$ for $m$ sufficiently large.
By sine law (\ref{sinelaw}), $l_2^{(m)}<l_1^{(m)}$.
\\

\noindent \textit{Case 2.} All of $\theta_1^{(m)}$, $\theta_2^{(m)}$
and $\theta_3^{(m)}$ converge to $+\infty$. In this case, since
$x_i^{(m)}+x_j^{(m)}=\theta_k^{(m)}$ converges to $+\infty$, by
taking a subsequence if necessary, at least two of $x_1^{(m)}$,
$x_2^{(m)}$ and $x_3^{(m)}$, say $x_1^{(m)}$ and $x_2^{(m)}$,
converge to $+\infty$. Therefore, $l_3$ is the only possible good
edge length, and $x_3^{(m)}$ converges to a finite number. For $m$
sufficiently large,
$\theta_1^{(m)}=x_2^{(m)}+x_3^{(m)}<x_1^{(m)}+x_2^{(m)}=\theta_3^{(m)}$
and
$\theta_2^{(m)}=x_1^{(m)}+x_3^{(m)}<x_1^{(m)}+x_2^{(m)}=\theta_3^{(m)}.$
By sine law (\ref{sinelaw}), $l_1^{(m)}<l_3^{(m)}$ and
$l_2^{(m)}<l_3^{(m)}$.\hfill{q.e.d}

\medskip
\noindent {\bf Lemma 3.4} \it
\begin{enumerate}[(a)]
\item If $\{(l_1,l_2,l_3)\}$ converges to $(+\infty,f_2,f_3)$ with $f_2,f_3\in\mathbb{R}$, then $\{(\theta_1,\theta_2,\theta_3)\}$ converges to
$(+\infty,0,0)$.
\item If $\{(l_1,l_2,l_3)\}$ converges to $(+\infty,+\infty,f_3)$ with $f_3\in\mathbb{R}$, then $\{\theta_3\}$ converges to
$0$.
\item If $\{(l_1,l_2,l_3)\}$ converges to $(+\infty,+\infty,+\infty)$, then we can take a subsequence such that at least two of $\{\theta_1\}$, $\{\theta_2\}$ and $\{\theta_3\}$ converge to $0$.
\end{enumerate}\rm
\medskip

\begin{figure}[htbp]\centering
\includegraphics[width=9cm]{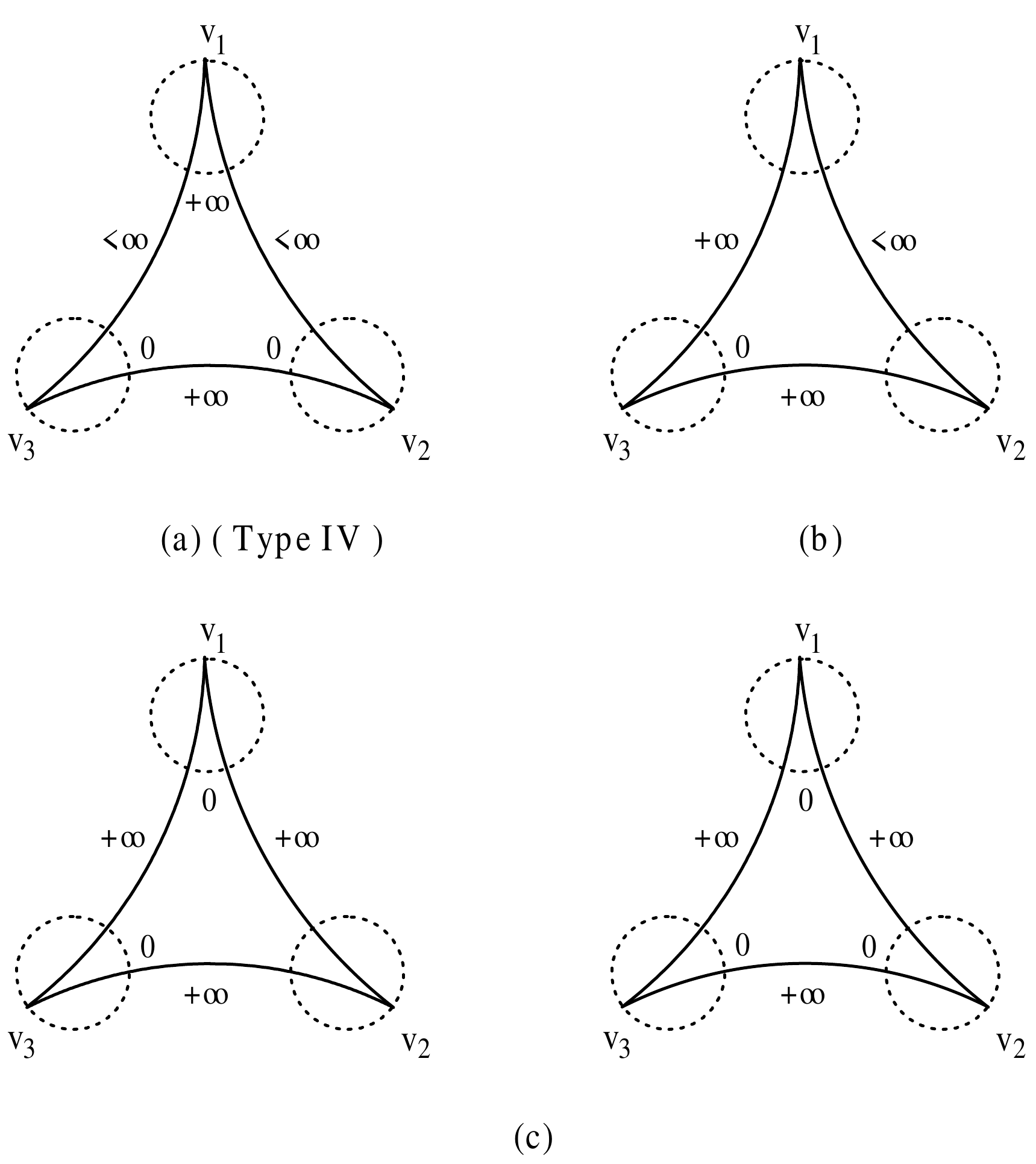}\\
\caption{Type IV and other types.}\label{type4}
\end{figure}

We call a converging sequence of decorated ideal triangles in (a) of
Lemma 3.4 a \emph{degenerated decorated ideal triangle of type IV}
(Figure \ref{type4}).
\\

\noindent \textit{Proof}: For (a), if $\{(l_1,l_2,l_3)\}$ converges
to $(+\infty,f_2,f_3)$, then by cosine law (\ref{cosinelaw}),
$\{\theta_1\}=\{e^{\frac{l_1-l_2-l_3}{2}}\}$ converges to $+\infty$,
$\{\theta_2\}=\{e^{\frac{l_2-l_1-l_3}{2}}\}$ converges to $0$, and
$\{\theta_3\}=\{e^{\frac{l_3-l_1-l_2}{2}}\}$ converges to $0$. For
(b), if $\{(l_1,l_2,l_3)\}$ converges to $(+\infty,+\infty,f_3)$,
then $\{\frac{l_3-l_1-l_2}{2}\}$ converges to $-\infty$, and
$\{\theta_3\}=\{e^{\frac{l_3-l_1-l_2}{2}}\}$ converges to $0$. For
(c), if $\{(l_1,l_2,l_3)\}$ converges to
$(+\infty,+\infty,+\infty)$, then we have by cosine law (\ref{cosinelaw})
that $\{\theta_1\theta_2\}=\{e^{-l_3}\}$ converges to $0$. Thus, by
taking a subsequence if necessary, at least one of $\{\theta_1\}$
and $\{\theta_2\}$, say $\{\theta_1\}$, converges to $0$. Since
$\{\theta_2\theta_3\}=\{e^{-l_1}\}$ converges to $0$ as well, by
taking a subsequence, at least one of $\{\theta_2\}$ and
$\{\theta_3\}$ converges to $0$. \hfill{q.e.d}


\bigskip
\centerline {\bf 4. The image of $\Psi_h$}
\bigskip

\noindent The image of $\Psi_h$ is described in Theorem 1.2. The
main task of this section is to give a proof of this theorem. To
show that the image of $\Psi_h$ is indeed $P_h(T)$, we make use of
the following propositions which are proved in this section.

\medskip
\noindent {\bf Proposition 4.1} \it
$\Psi_h(T_c(S)\times\mathbb{R}_{>0}^V)\subset P_h(T)$ for all
$h\in\mathbb{R}$.\rm
\medskip

\noindent {\bf Proposition 4.2} \it For all $h\in\mathbb{R}$, the
image $\Psi_h(T_c(S)\times\mathbb{R}_{>0}^V)$ is closed in $P_h(T)$
.\rm
\medskip

\noindent \textit{Proof of Theorem 1.2}: Let
$P(T)$ be defined as in Theorem 1.2. For $h\geqslant0$, $P(T)=P_h(T)$ is determined by finitely many strict linear inequalities corresponding to the fundamental edge loops and hence is an open convex polytope independent of $h$. For $h<0$, $P_h(T)$ is likewise determined by fundamental edge loops and fundamental edge paths. Moreover, since each
edge $e$ can be regarded as a fundamental edge path, conditions (a)
and (b) imply that
$-2\int_0^{+\infty}e^{ht^2}dt<z(e)<2\int_0^{+\infty}e^{ht^2}dt$ for
each $e\in E$. Thus, $P_h(T)$ is bounded. The monotonicity of the
function $f(h)=\int_0^{+\infty}e^{ht^2}dt$ implies that
$P_h(T)\subset P_{h'}(T)$ if $h<h'$, and the fact that
$\lim_{h\rightarrow-\infty}f(h)=\lim_{h\rightarrow-\infty}\sqrt{\frac{\pi}{-2h}}=0$
implies that $\bigcap_{h\in\mathbb{R}_{<0}}P_h(T)=\emptyset$. By
Theorem 1.1 and the Invariance of Domain Theorem,
$\Psi_h(T_c(S)\times\mathbb{R}_{>0}^V)$ is open in $P_h(T)$. By
Proposition 4.2, $\Psi_h(T_c(S)\times\mathbb{R}_{>0}^V)$ is closed
in $P_h(T)$. Connectedness of $P_h(T)$ therefore implies that
$\Psi_h(T_c(S)\times\mathbb{R}_{>0}^V)= P_h(T)$.\hfill{q.e.d}
\\

The following Lemma 4.3 will be used in the proof of Propositions
4.1 and 4.2.

\medskip
\noindent {\bf Lemma 4.3} \it If $r\in \mathbb{R}$ and $x>0$, then

\begin{enumerate}[(a)]
\item for each $h\in\mathbb{R}$,
$$\int_0^{x+r}e^{ht^2}dt+\int_0^{x-r}e^{ht^2}dt>0,$$
\item for each $h\geqslant0$,$$\int_0^{x+r}e^{ht^2}dt+\int_0^{x-r}e^{ht^2}dt\geqslant2\int_0^xe^{ht^2}dt.$$
\end{enumerate}\rm
\medskip

\noindent \textit{Proof}: For (a), let $f(x)=\int_0^{x+r}e^{ht^2}dt+\int_0^{x-r}e^{ht^2}dt$. Since $f'(x)=e^{h(x+r)^2}+e^{h(x-r)^2}>0$, the function $f$ is strictly increasing, hence $f(x)>f(0)=0$ for $x>0$. For (b), let
$g(x)=\int_0^{x+r}e^{ht^2}dt+\int_0^{x-r}e^{ht^2}dt-2\int_0^xe^{ht^2}dt.$
We have that $g(0)=0$ and
$g'(x)=e^{h(x+r)^2}+e^{h(x-r)^2}-2e^{hx^2}\geqslant0$. The last
inequality follows from the convexity of the function $F(t)=e^{ht^2}$ for
$h\geqslant0$. Since $g$ is increasing, $g(x)\geqslant g(0)=0$ for
$x>0$. \hfill{q.e.d}
\\

\noindent \textit{Proof of Proposition 4.1}: For $h\geqslant0$, fix
a decorated hyperbolic metric $(d,r)\in
T_c(S)\times\mathbb{R}_{>0}^V$. For any fundamental edge loop
$(e_1,t_1,\dots,e_k,t_k)$, let $a_i$ be the generalized angle adjacent
to $e_i$ and $e_{i+1}$ (where $e_{k+1}=e_1$). Let the generalized
angles of $t_i$ facing $e_i$ and $e_{i+1}$ respectively be $b_i$ and $c_i$. By definition, the contribution of $\sum_{i=1}^kz(e_i)$ from $t_i$
is

$$\int_0^{\frac{a_i+b_i-c_i}{2}}e^{ht^2}dt+\int_0^{\frac{a_i+c_i-b_i}{2}}e^{ht^2}dt,$$ which is strictly larger than $0$ from Lemma 4.3\,(a)
since $a_i>0$.
\\

For $h<0$, let $e$ be any edge in the ideal triangulation $T$, and
let $a$ and $a'$ be the generalized angles facing $e$. Let $b$, $c$,
$b'$ and $c'$ be the generalized angles adjacent to $e$. Then

$$\Psi_h(d,r)(e)=\int_0^{\frac{b+c-a}{2}}e^{ht^2}dt+\int_0^{\frac{b'+c'-a'}{2}}e^{ht^2}dt<2\int_0^{+\infty}e^{ht^2}dt.$$
Thus, condition (a) in the definition of $P_h(T)$ is satisfied.
Given a fundamental edge path $(t_0,e_0,t_1,\dots,e_n,t_n)$, let
$\theta_i$ be the generalized angle in $t_i$ adjacent to $e_i$ and
$e_{i+1}$ for $i=1,\dots,n-1$, and let $\beta_i$ and $\gamma_i$ respectively be the
generalized angles of $t_i$ facing $e_i$ and $e_{i+1}$. Denote by
$a_0$ the generalized angle of $t_0$ facing $e_0$, and by $a_n$ the
generalized angle of $t_n$ facing $e_n$. Let $b_0$ and $c_0$ be the
generalized angles of $t_0$ adjacent to $e_0$, and let $b_n$ and $c_n$ be the
generalized angles of $t_n$ adjacent to $e_n$. We have

\begin{equation*}
\begin{split}
&\sum_{i=1}^n\Psi_h(d,r)(e_i)\\
=&\int_0^{\frac{b_0+c_0-a_0}{2}}e^{ht^2}dt+\sum_{i=1}^{n-1}\big(\int_0^{\frac{\theta_{i}+\gamma_{i}-\beta_{i}}{2}}e^{ht^2}dt+\int_0^{\frac{\theta_i+\beta_i-\gamma_i}{2}}e^{ht^2}dt\big)+\int_0^{\frac{b_n+c_n-a_n}{2}}e^{ht^2}dt\\
>&\int_0^{\frac{b_0+c_0-a_0}{2}}e^{ht^2}dt+\int_0^{\frac{b_n+c_n-a_n}{2}}e^{ht^2}dt\\
>&-2\int_0^{+\infty}e^{ht^2}dt,
\end{split}
\end{equation*}
where the first inequality is by Lemma 4.3\,(a). Thus, condition (b)
is satisfied. Given a fundamental edge loop $(e_1,t_1,\dots,e_n,t_n)$
with $e_{n+1}=e_1$, let $\theta_i$ for $i=1,\dots,n$ be the
generalized angle in $t_i$ adjacent to $e_i$ and $e_{i+1}$, and let
$\beta_i$ (resp. $\gamma_i$) be the generalized angle in $t_i$
facing $e_i$ (resp. $e_{i+1}$). Again by Lemma 4.3\,(a),

\begin{equation*}
\sum_{i=1}^n\Psi_h(d,r)(e_i)=\sum_{i=1}^{n}\big(\int_0^{\frac{\theta_{i}+\gamma_{i}-\beta_{i}}{2}}e^{ht^2}dt+\int_0^{\frac{\theta_i+\beta_i-\gamma_i}{2}}e^{ht^2}dt\big)>0.
\end{equation*}
Thus, condition (c) is satisfied, and
$\Psi_h(T_c(S)\times\mathbb{R}_{>0}^V)\subset
P_h(T)$.\hfill{q.e.d}
\\

To prove Proposition 4.2, we use Penner's length parametrization.
For each sequence $\{l^{(m)}\}$ in $\mathbb{R}^E$ such that
$\{\Psi_h(l^{(m)})\}$ converges to a point $z\in P(T)$, we claim
that $\{l^{(m)}\}$ contains a subsequence converging to a point in
$\mathbb{R}^E$. Let $\theta^{(m)}$ be the generalized angles of the
decorated ideal triangles in $(S,T)$ in the decorated hyperbolic
metric $l^{(m)}$. By taking a subsequence if necessary, we may
assume that $\{l^{(m)}\}$ converges in $[-\infty,+\infty]^E$ and
that for each generalized angle $\theta_i$, the limit
$\lim_{m\rightarrow\infty}\theta_i^{(m)}$ exists in $[0,+\infty]$.
In the case that $h\geqslant0$, we need the following

\medskip
\noindent {\bf Lemma 4.4} \it If $h\geqslant0$, then
$\lim_{m\rightarrow\infty}\theta_i^{(m)}\in[0,+\infty)$ for all
$i$.\rm
\medskip

\noindent \textit{Proof}: Suppose to the contrary that
$\lim_{m\rightarrow\infty}\theta_1^{(m)}=+\infty$ for some
generalized angle $\theta_1$. Let $e_2$ and $e_3$ be the edges
adjacent to $\theta_1$ in the triangle $t_1$, and $\theta_2$ and
$\theta_3$ respectively be the generalized angles facing $e_2$ and $e_3$. Take a
fundamental edge loop $(e_{n_1},t_{n_1},\dots,e_{n_k},t_{n_k})$
containing $(e_2,t_1,e_3)$. By Lemma 4.3, we have

\begin{equation*}
\begin{split}
\sum_{i=1}^kz(e_{n_i})=&\lim_{m\rightarrow\infty}\sum_{i=1}^k\Psi_h(l^{(m)})(e_{n_i})\\
\geqslant&\lim_{m\rightarrow\infty}\big(\int_0^{\frac{\theta_1^{(m)}+\theta_2^{(m)}-\theta_3^{(m)}}{2}}e^{ht^2}dt+\int_0^{\frac{\theta_1^{(m)}+\theta_3^{(m)}-\theta_2^{(m)}}{2}}e^{ht^2}dt\big)\\
\geqslant&\lim_{m\rightarrow\infty}2\int_0^{\frac{\theta_1^{(m)}}{2}}e^{ht^2}dt\\
=&+\infty.
\end{split}
\end{equation*}
This contradicts the assumption that $z\in P(T)$. \hfill{q.e.d}
\\

\noindent \textit{Proof of Proposition 4.2}: For $h\geqslant0$, by
taking a subsequence of $\{l^{(m)}\}$, we may assume that
$\lim_{m\rightarrow\infty}l^{(m)}=l\in[-\infty,+\infty]^E$. If $l$
were not in $\mathbb{R}^E$, then there would exist an edge $e\in E$ so
that $l(e)=\pm\infty$. Let $\Delta$ be a decorated ideal triangle
adjacent to $e$, and let $\theta_1^{(m)}$ and $\theta_2^{(m)}$ be
the generalized angles in $\Delta$ adjacent to $e$ in the metric
$l^{(m)}$. By (\ref{cosinelaw}),

\begin{equation*}
e^{l^{(m)}(e)}=\frac{1}{\theta_1^{(m)}\theta_2^{(m)}},
\end{equation*}

\noindent and $\theta_i^{(m)}\in(0,+\infty)$ for $i=1,2$.
\\

\noindent \textit{Case 1} If $l(e)=-\infty$, then $e^{l(e)}=0$. By
the identity above, one of $\lim_{m\rightarrow\infty}\theta_i^{(m)}$
for $i=1,2$ must be $+\infty$. This contradicts Lemma 4.4.
\\

\noindent \textit{Case 2} If $l(e)=+\infty$, then
$e^{l(e)}=+\infty$. By the identity above, one of
$\lim_{m\rightarrow\infty}\theta_i^{(m)}$ for $i=1,2$ must be zero.
 Suppose without loss of generality that $\lim_{m\rightarrow\infty}\theta_1^{(m)}=0$. Let $e_1$ be the
edge in the decorated ideal triangle $\Delta$ opposite to
$\theta_2$, and let $\theta_3$ be the generalized angle in $\Delta$
facing $e$. By (\ref{cosinelaw}), we have
\begin{equation*}
e^{l^{(m)}(e_1)}=\frac{1}{\theta_1^{(m)}\theta_3^{(m)}}.
\end{equation*}
By Lemma 4.4, $\theta_3^{(m)}$ is bounded above, hence
$l(e_1)=+\infty$. For any decorated ideal triangle $\Delta$ adjacent
to $e$ with $l(e)=+\infty$, we have an edge $e_1$ in $\Delta$ and a
generalized angle $\theta_1$ adjacent to $e$ and $e_1$ so that
$l(e_1)=+\infty$ and $\lim_{m\rightarrow\infty}\theta_1^{(m)}=0.$
Applying this logic to $e_1$ and the decorated ideal triangle
$\Delta_1$ adjacent to $e_1$ other than $\Delta$, we obtain the next
angle $\theta_2$ and edge $e_2$ in $\Delta_1$ so that
$l(e_2)=+\infty$ and $\lim_{m\rightarrow\infty}\theta_2^{(m)}=0$.
Since there are only finitely many edges and triangles, this yields
a fundamental edge loop $(e_k,\Delta_k,\dots,e_n,\Delta_n)$ in $T$
such that $l(e_i)=+\infty$ for $i=k,\dots,n$ and
$\lim_{m\rightarrow\infty}\theta_i^{(m)}=0$, where $\theta_i$ is the
generalized angle in $\Delta_{i-1}$ adjacent to $e_{i-1}$ and $e_i$.
Denote respectively by $\beta_i$ and $\gamma_i$ the generalized angles of
$\Delta_{i-1}$ facing $e_{i-1}$ and $e_i$, and let
$\bar{\beta}_i=\lim_{m\rightarrow\infty}\beta_i^{(m)}$ and
$\bar{\gamma}_i=\lim_{m\rightarrow\infty}\gamma_i^{(m)}$. By Lemma
4.4, both $\bar{\beta}_i$ and $\bar{\gamma}_i$ are finite real
numbers, and we have
\begin{equation*}
\begin{split}
\sum_{i=k}^nz(e_i)=&\lim_{m\rightarrow\infty}\sum_{i=k}^n\Psi_h(l^{(m)})(e_i)\\
=&\lim_{m\rightarrow\infty}\sum_{i=k}^n\big(\int_0^{\frac{\theta_i^{(m)}+\beta_i^{(m)}-\gamma_i^{(m)}}{2}}e^{ht^2}dt+\int_0^{\frac{\theta_i^{(m)}+\gamma_i^{(m)}-\beta_i^{(m)}}{2}}e^{ht^2}dt\big)\\
=&\sum_{i=k}^n\big(\int_0^{\frac{\bar{\beta}_i-\bar{\gamma}_i}{2}}e^{ht^2}dt+\int_0^{\frac{\bar{\gamma}_i-\bar{\beta}_i}{2}}e^{ht^2}dt\big)\\
=&\,0.
\end{split}
\end{equation*}
This contradicts the assumption that $z\in P(T)$.
\\

For $h<0$ and each sequence $\{l^{(m)}\}$ in $\mathbb{R}^E$ so that
$\{\Psi_h(l^{(m)})\}$ converges to a point $z\in P_h(T)$, we claim
that $\{l^{(m)}\}$ contains a subsequence converging to a point in
$\mathbb{R}^E$. By taking a subsequence if necessary, we may assume
that $\{l^{(m)}\}$ converges to $l\in[-\infty,+\infty]^E$. If $l$
were not in $\mathbb{R}^E$, there would exist an edge $e$ so that
$l(e)=\pm\infty$.
\\

\noindent \textit{Case 1.} If $l(e)=-\infty$ for some $e\in E$, then
there is a degenerated decorated ideal triangle $\Delta$ of type I,
II or III. By Corollary 3.2, there is a bad edge $e_1$ in $\Delta$.
Let $\Delta_1$ be the other decorated ideal triangle adjacent to
$e_1$, and let $x_0$ and $x_1$ respectively be the $x$-invariants of $e_1$ in
$\Delta$ and $\Delta_1$. If $e_1$ is bad in $\Delta_1$, then

\begin{equation*}
z(e_1)=\lim_{m\rightarrow\infty}\Psi_h(l^{(m)})(e_1)
=\lim_{m\rightarrow\infty}\big(\int_0^{x_0^{(m)}}e^{ht^2}dt+\int_0^{x_1^{(m)}}e^{ht^2}dt\big)
=2\int_0^{+\infty}e^{ht^2}dt,
\end{equation*}
which contradicts the assumption that $z\in P_h(T)$. Therefore $e_1$
has to be a good edge in $\Delta_1$. Since $l(e_1)=-\infty$, the
decorated triangle $\Delta_1$ is degenerated of type I, II or III.
By Corollary 3.2, there is a bad edge $e_2$ in $\Delta_1$. For the
same reason, $e_2$ has to be good in the other decorated ideal
triangle $\Delta_2$ adjacent to $e_2$, and there is a bad edge $e_3$
in $\Delta_2$. Serially applying this logic and using that there are finitely many edges, we find an edge loop
$(e_k,\Delta_k,\dots,e_n,\Delta_n)$ with $e_{n+1}=e_k$ so that for
each $i=k,\dots,n$ the edge $e_i$ is good in $\Delta_i$ and the edge
$e_{i+1}$ is bad in $\Delta_i$. By Lemma 3.3, we can take a
subsequence so that $l^{(m)}(e_i)>l^{(m)}(e_{i+1})$ for $m$
sufficiently large. Thus, we have $l^{(m)}(e_k)>l^{(m)}(e_{n+1})$,
which contradicts that $e_{n+1}=e_k$.
\\

In light of Case 1, we may assume that $l\in (-\infty,+\infty]^E$.
\\

\noindent \textit{Case 2.} If $l(e)=+\infty$ for some $e\in E$, let
$\Delta_1$ be a decorated ideal triangle adjacent to $e$. If
$\Delta_1$ is not of type IV, then by Lemma 3.4, there is an edge
$e_1$ of $\Delta_1$ and an generalized angle $\theta_1$ adjacent to
$e$ and $e_1$ so that $l(e_1)=+\infty$ and
$\lim_{m\rightarrow\infty}\theta_1^{(m)}=0$ (see Figure
\ref{type4}). The other decorated ideal triangle $\Delta_2$ adjacent to $e_1$ is either of type IV or contains an edge
$e_2$ and a generalized angle $\theta_2$ adjacent to $e_1$ and $e_2$
so that $l(e_2)=+\infty$ and
$\lim_{m\rightarrow\infty}\theta_2^{(m)}=0$. Again, the serial application of this procedure terminates with an edge $e_p$ and a decorated
ideal triangle $\Delta_{p+1}$ adjacent to $e_p$ so that
$l(e_p)=+\infty$ and $\Delta_{p+1}$ is of type IV, or since there
are only finitely many edges, produces a fundamental edge
loop $(e_k,\Delta_k,\dots,e_n,\Delta_n)$ such that $l(e_i)=+\infty$
for $i=k,\dots,n$ and $\lim_{m\rightarrow\infty}\theta_i^{(m)}=0$,
where $\theta_i$ is the generalized angle in $\Delta_i$ adjacent to
$e_i$ and $e_{i+1}$. If it yields such a fundamental edge loop
$(e_k,\Delta_k,\dots,e_n,\Delta_n)$, denote by $\beta_i$ (resp.
$\gamma_i$) the generalized angle in $\Delta_i$ facing $e_i$ (resp.
$e_{i+1}$) for $i=k,\dots,n$. Let
$\bar{\beta}_i=\lim_{m\rightarrow\infty}\beta_i^{(m)}$ and
$\bar{\gamma}_i=\lim_{m\rightarrow\infty}\gamma_i^{(m)}$, so that
\begin{equation*}
\begin{split}
\sum_{i=k}^nz(e_i)=&\lim_{m\rightarrow\infty}\sum_{i=1}^k\Psi_h(l^{(m)})(e_i)\\
=&\lim_{m\rightarrow\infty}\sum_{i=1}^k\big(\int_0^{\frac{\theta_i^{(m)}+\beta_i^{(m)}-\gamma_i^{(m)}}{2}}e^{ht^2}dt+\int_0^{\frac{\theta_i^{(m)}+\gamma_i^{(m)}-\beta_i^{(m)}}{2}}e^{ht^2}dt\big)\\
=&\sum_{i=1}^k\big(\int_0^{\frac{\bar{\beta}_i-\bar{\gamma}_i}{2}}e^{ht^2}dt+\int_0^{\frac{\bar{\gamma}_i-\bar{\beta}_i}{2}}e^{ht^2}dt\big)\\
=&\,0,
\end{split}
\end{equation*}
which contradicts the assumption that $z\in P_h(T)$. If it terminates with $e_p$ and $\Delta_{p+1}$ of type IV, then we consider the other decorated
ideal triangle $\Delta_0$ adjacent to $e$. If $\Delta_0$ is not of
type IV, then it contains an edge $e_{-1}$ and a generalized angle
$\theta_0$ adjacent to $e_{-1}$ and $e$ so that $l(e_{-1})=+\infty$
and $\lim_{m\rightarrow\infty}\theta_0^{(m)}=0$. As before, either there is a fundamental edge loop, contradicting the assumption that $z\in P_h(T)$, or the procedure terminates with an edge $e_{-q}$ and a decorated ideal triangle
$\Delta_{-q}$ adjacent to $e_{-q}$ so that $l(e_{-q})=+\infty$ and
$\Delta_{-q}$ is of type IV. If the procedure stops at $e_{-q}$ and
$\Delta_{-q}$ of type IV, we get a fundamental edge path
$(\Delta_{-q},e_{-q},\dots,e_{p},\Delta_{p+1})$, where $e_0=e$, such
that $\Delta_{-q}$ and $\Delta_p$ are of type IV with
$l(e_{-q})=+\infty$ and $l(e_p)=+\infty$, and
$\lim_{m\rightarrow\infty}\theta_i^{(m)}=0$, where $\theta_i$ is the
generalized angle of $\Delta_i$ adjacent to $e_{i-1}$ and $e_i$ for
$i=1-q,\dots,p$. Denote by $a_{-q}$ the generalized angle of
$\Delta_{-q}$ facing $e_{-q}$, and by $a_p$ the generalized angle of
$\Delta_{p+1}$ facing $e_p$. Let $b_{-q}$ and $c_{-q}$ be the
generalized angles of $\Delta_{-q}$ adjacent to $e_{-q}$, and let $b_p$ and $c_p$ be
the generalized angles of $\Delta_{p+1}$ adjacent to $e_p$. We find
\begin{equation*}
\begin{split}
\sum_{i=-q}^pz(e_i)=&\lim_{m\rightarrow\infty}\sum_{i=-q}^{p}\Psi_h(l^{(m)})(e_i)\\
=&\lim_{m\rightarrow\infty}\big(\int_0^{\frac{b_{-q}^{(m)}+c_{-q}^{(m)}-a_{-q}^{(m)}}{2}}e^{ht^2}dt+\int_0^{\frac{b_p^{(m)}+c_p^{(m)}-a_p^{(m)}}{2}}e^{ht^2}dt\\
&+\sum_{i=1-q}^{p}\big(\int_0^{\frac{\theta_i^{(m)}+\beta_i^{(m)}-\gamma_i^{(m)}}{2}}e^{ht^2}dt+\int_0^{\frac{\theta_i^{(m)}+\gamma_i^{(m)}-\beta_i^{(m)}}{2}}e^{ht^2}dt\big)\big)\\
=&\int_0^{-\infty}e^{ht^2}dt+\int_0^{-\infty}e^{ht^2}dt+\sum_{i=1-q}^p\big(\int_0^{\frac{\bar{\beta}_i-\bar{\gamma}_i}{2}}e^{ht^2}dt+\int_0^{\frac{\bar{\gamma}_i-\bar{\beta}_i}{2}}e^{ht^2}dt\big)\\
=&-2\int_0^{+\infty}e^{ht^2}dt,
\end{split}
\end{equation*}
which contradicts the assumption that $z\in P_h(T)$. \hfill{q.e.d}


\bigskip
\centerline {\bf 5. Uniqueness of the energy function}
\bigskip

\noindent Let $\Delta$ be a decorated ideal triangle with edge
lengths $l_1$, $l_2$, $l_3$ with opposite generalized angles
$\theta_1$, $\theta_2$, $\theta_3$ and set
$x_i=\frac{\theta_j+\theta_k-\theta_i}{2}$ for $i,j,k=1,2,3$. The
following theorem shows that $\Psi_h$ is the unique possible
deformation of Penner's simplicial coordinate by using the variational
principle stated in Section 2.

\medskip
\noindent {\bf Theorem 5.1} \it Let $\mu$ and $u$ be two
non-constant smooth functions. Up to an overall scale, there is a unique closed
$1$-form $\omega=\sum_{i=1}^3\mu(x_i)du(l_i)$ which is given by
$$w_h=\sum_{i=1}^3\int^{x_i}e^{ht^2}dtd\big(\int^{l_i}e^{-he^{-t}}dt\big)$$
for some $h\in \mathbb{R}$.\rm
\medskip

The proof of Theorem 5.1 makes use of the following lemma.

\medskip
\noindent {\bf Lemma 5.2} \it Let $f$ and $g$ be
two non-constant smooth functions on $\mathbb{R}$. If
$\frac{f(x_i)}{g(l_j)}$ is symmetric in $i,j=1,2$, then there are constants $h$, $c_1$ and $c_2$ so that
$$f(t)=e^{ht^2+c_1}\quad \text{and}\quad g(t)=e^{-he^{-t}+c_2}.$$\rm
\medskip

\noindent \textit{Proof}: By taking $\frac{\partial}{\partial l_k}$
in the equality $\frac{f(x_i)}{g(l_j)}=\frac{f(x_j)}{g(l_i)}$, we
have $\frac{f'(x_i)}{g(l_j)}\frac{\partial x_i}{\partial
l_k}=\frac{f'(x_j)}{g(l_i)}\frac{\partial x_j}{\partial l_k}$ for $i,j,k=1,2,3$. We
deduce from (\ref{cosinelaw}) that $\frac{\partial x_i}{\partial
l_j}=\frac{x_k}{2}$, so
$\frac{f'(x_i)}{g(l_j)}\frac{x_j}{2}=\frac{f'(x_j)}{g(l_i)}\frac{x_i}{2}.$
Thus,
$\frac{f'(x_i)}{f'(x_j)}\frac{x_j}{x_i}=\frac{g(l_j)}{g(l_i)}=\frac{f(x_i)}{f(x_j)},$
which implies
$\frac{f'(x_i)}{f(x_i)}\frac{1}{x_i}=\frac{f'(x_j)}{f(x_j)}\frac{1}{x_j}$ and $\frac{f'(t)}{f(t)}\frac{1}{t}=2h_1$ for some
$h_1\in\mathbb{R}$. Solving this ordinary differential equation for
$f$, we find
\begin{equation*}
f(t)=e^{h_1t^2+c_1}
\end{equation*}
for some $c_1\in\mathbb{R}$.
By taking $\frac{\partial}{\partial x_k}$ in the equality
$\frac{g(l_i)}{f(x_j)}=\frac{g(l_j)}{f(x_i)}$, we have
$\frac{g'(l_i)}{f(x_j)}\frac{\partial l_i}{\partial
x_k}=\frac{g'(l_j)}{f(x_i)}\frac{\partial l_j}{\partial x_k}$ for $i,j,k=1,2,3$. From
(\ref{cosinelaw}) again, we deduce that $\frac{\partial l_i}{\partial x_j}=-\frac{1}{\theta_k}$, so
$-\frac{g'(l_i)}{f(x_j)}\frac{1}{\theta_j}=-\frac{g'(l_j)}{f(x_i)}\frac{1}{\theta_i}$.
Thus, $\frac{g'(l_i)}{g'(l_j)}\frac{e^{l_i}}{e^{l_j}}=\frac{g'(l_i)}{g'(l_j)}\frac{\theta_i}{\theta_j}=\frac{f(x_j)}{f(x_i)}=\frac{g(l_i)}{g(l_j)},$
which implies
$\frac{g'(l_i)}{g(l_i)}e^{l_i}=\frac{g'(l_j)}{g(l_j)}e^{l_j}$ and
$\frac{g'(t)}{g(t)}e^{t}=h_2$ for some} $h_2\in\mathbb{R}$.
Solving this ordinary differential equation for $g$, we find
\begin{equation*}
g(t)=e^{-h_2e^{-t}+c_2}
\end{equation*}
 for some $c_1\in\mathbb{R}$. From $f(t)=e^{h_1t^2+c_1}$ and the equality
$\frac{f(x_i)}{g(l_j)}=\frac{f(x_j)}{g(l_i)}$, we conclude that
$h_1=h_2$. \hfill{q.e.d}
\\

\noindent \textit{Proof of Theorem 5.1}: The differential $1$-form
$\omega=\sum_{i=1}^3\mu(x_i)du(l_i)$ is closed if and only if
$\frac{\partial \mu(x_i)}{\partial
u(l_j)}=\frac{\mu'(x_i)}{u'(l_j)}\frac{\partial x_i}{\partial l_j}$
is symmetric in $i$ and $j$. Since $\frac{\partial x_i}{\partial
l_j}=\frac{\partial x_j}{\partial l_i}=\frac{x_k}{2}$, $\omega$ is closed if and only if
$\frac{\mu'(x_i)}{u'(l_j)}$ is symmetric in $i$ and $j$. By Lemma
5.2, if $\frac{\mu'(x_i)}{u'(l_j)}$ is symmetric in $i$ and $j$,
then $\mu'(x_i)=e^{hx_i^2+c_1}$ and $u'(l_i)=e^{-he^{-l_i}+c_2}$ for
some constants $h$, $c_1$ and $c_2$. \hfill{q.e.d}


\bigskip
\centerline{\bf 6. $\Psi_h$ and the Delaunay decomposition}
\bigskip

\noindent We first review the construction of the Delaunay
decomposition associated to a decorated hyperbolic metric following
Bowditch-Epstein\,[BE]. Suppose $S$ is a punctured surface with a
set of ideal vertices $V$, and let $(d,r)$ be a decorated hyperbolic
metric on $S$ so that the horodisks associated to the ideal vertices
do not intersect. Let $B_v$ be the horodisks associated to the ideal
vertex $v$, and let $B=\bigcup_{v\in V}B_v$. The \emph{spine}
$\Gamma_{d,r}$ of $S$ is the set of points in $S$ which have at
least two distinct shortest geodesics to $\partial B$. The spine
$\Gamma_{d,r}$ is shown\,[BE] to be a graph whose edges are geodesic
arcs on $S$.
\\

Let $e_1^*,\dots,e_N^*$ be the edges of $\Gamma_{d,r}$. By
construction each interior point of an edge $e_i^*$ has exactly two
distinct shortest geodesics to $\partial B$. For each edge $e_i^*$,
there are two horodisks $B_1$ and $B_2$ (possibly coincide) so that
points in the interior of $e_i^*$ have precisely two shortest
geodesics to $\partial B_1$ and $\partial B_2$. Let $e_i$ be the
shortest geodesic from $\partial B_1$ to $\partial B_2$. It is known
that $e_i$ intersects $e_i^*$ perpendicularly, and $\{e_1,\dots,e_N\}$
are disjoint. The components of $S\setminus\{e_1,\dots,e_N\}$ consists
of decorated polygons (ideal polygons with horodisks associated to
the ideal vertices) which are the $2$-cells of the \emph{Delaunay
decomposition} $\Sigma_{d,r}$. The $1$-cells of $\Sigma_{d,r}$
consist of the edges $\{e_1,\dots,e_N\}$ and the arcs on $\partial B$
which are the intersection of $\partial B$ with the ideal polygons.
For a generic decorated hyperbolic metric $(d,r)$, each $2$-cell of
$\Sigma_{d,r}$ is a decorated ideal triangle, and $\Sigma_{d,r}$ is
a decorated ideal triangulation of $S$.
\\

Let $D$ be a $2$-cell of $\Sigma_{d,r}$. We call the hyperbolic
circle on $S$ tangent to all arcs of $D\cap
\partial B$ the \emph{inscribed circle} of $D$. By the construction
of the Delaunay decomposition, for each $2$-cell $D$ of
$\Sigma_{d,r}$, there is exactly one vertex $v^*$ of the spine
$\Gamma_{d,r}$ lying in the interior of $D$. Moreover, $v^*$ is of
equal distance to all arcs of $D\cap \partial B$, hence is the
center of the inscribed circle of $D$. Thus, the center of the
inscribed circle of each $2$-cell $D$ of the Delaunay decomposition
is in the interior of $D$. We need the following proposition of Penner\,[P3] whose proof is included here to the convenience of the readers.

\medskip
\noindent {\bf Lemma 6.1}\,\bf([P3]) \it Suppose $\Delta$ is a decorated ideal
triangle with edge lengths $l_i>0$ and opposite generalized angles
$\theta_i$ for $i=1,2,3$. Then
$x_i=\frac{\theta_j+\theta_k-\theta_i}{2}>0$ for $i=1,2,3$ if and
only if the center of the inscribed circle of $\Delta$ is in the
interior of $\Delta$.\rm
\medskip

\begin{figure}[htbp]\centering
\includegraphics[width=10cm]{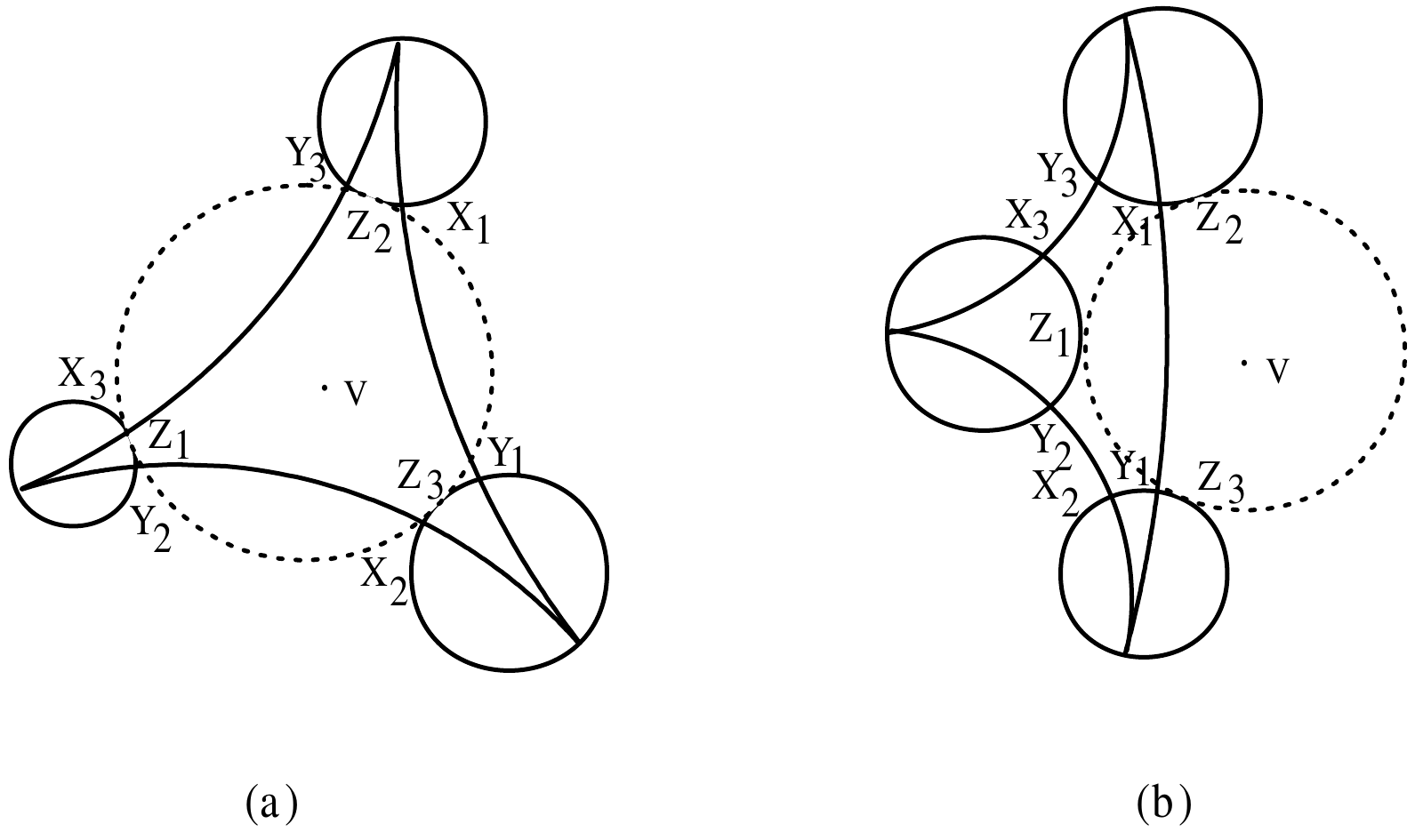}\\
\caption{The inscribed circle.}\label{figure2}
\end{figure}

\noindent \textit{Proof}: For $i=1,2,3$ let $B_i$ be the horodisks
associated to the ideal vertices of $\Delta$, and let $Z_i$ be the
point of tangency of the inscribe circle of $\Delta$ and $\partial
B_i$. Label the intersection of the horodisks and the edges of
$\Delta$ by $X_1,Y_1,X_2,Y_2,X_3$ and $Y_3$ cyclically as in Figure
\ref{figure2}(a). For two points $A$ and $B$ in the hyperbolic plane
$\mathbb{H}^2$, let $AB$ be the geodesic segment connecting $A$ and
$B$, and let $|AB|$ the length of $AB$. If the center $v$ of the
inscribed circle is in the interior of $\Delta$, then
$x_i=|X_iZ_{i+1}|>0$ for $i=1,2,3$. If $v$ is on $X_iY_i$, or $v$
and $\Delta$ are on different sides of $X_iY_i$ for some
$i\in\{1,2,3\}$, then $x_i=-|X_iZ_{i+1}|\leqslant0$. See Figure
\ref{figure2}\,(b). \hfill{q.e.d}
\\

\noindent \textit{Proof of Theorem 1.3}: Let $(d,r)$ be a decorated
hyperbolic metric so that the associated Delaunay decomposition
$\Sigma_{d,r}$ is a decorated ideal triangulation of $S$. For each
edge $e$ of $\Sigma_{d,r}$, let $\Delta$ and $\Delta'$ be the
decorated ideal triangles adjacent to $e$, and let $\theta_1$ and
$\theta_1'$ respectively be the generalized angles of $\Delta$ and $\Delta'$
facing $e$, and $\theta_2$, $\theta_3$, $\theta_2'$ and $\theta_3'$
be the generalized angles adjacent to $e$. Let
$x(e)=\frac{\theta_2+\theta_3-\theta_1}{2}$ and
$x'(e)=\frac{\theta_2'+\theta_3'-\theta_1'}{2}$. From Lemma 6.1
and the fact that the center of the inscribed circle of each
$2$-cell of the Delaunay decomposition is in the interior of the
$2$-cell, we conclude that $x(e)$ and $x'(e)$ are positive, and
$$\Psi_h(d,r)(e)=\int_0^{x(e)}e^{ht^2}dt+\int_0^{x'(e)}e^{ht^2}dt>0.$$
\\

On the other hand, if $T$ is an ideal triangulation of $S$ such that
$\Psi_h(d,r)(e)\leqslant0$ for some edge $e$, then at least one of
$x(e)$ and $x'(e)$, say $x(e)$, is less than or equal to zero. By
Lemma 6.1, the center of the inscribed circle of $\Delta$ is not in
the interior of $\Delta$. Since the center of the inscribed circle
of each $2$-cell of the Delaunay decomposition has to be in the
interior of the $2$-cell, $T$ cannot be the Delaunay decomposition
$\Sigma_{d,r}$ of $S$. \hfill{q.e.d}


\bigskip
\centerline{\bf 7. Further questions}
\bigskip

\noindent 1. Suppose $\Delta$ is a decorated ideal triangle with
edge lengths $l_1$, $l_2$ and $l_3$ and opposite generalized angles
$\theta_1$, $\theta_2$ and $\theta_3$. For each $h\neq-1$, the
differential $1$-form
$\omega_h=\sum_{i=1}^3\theta_i^{h+1}de^{-(h+1)l_i}$ is closed in
$\mathbb{R}^3$. However, the primitive $F_h(u)=\int_0^u\omega_h$
is not strictly concave on $\mathbb{R}^3$. Let $(S,T)$ be an ideally
triangulated punctured surface. For each $h\neq-1$, we define a
map $\Phi_h\co T_c(S)\times\mathbb{R}_{>0}^V\rightarrow\mathbb{R}^E$
by

$$\Phi_h(d,r)(e)=\theta^{h+1}+\theta'^{h+1},$$
where $\theta$ and $\theta'$ are the generalized angles facing $e$.
To the best of the author's knowledge, there is no counterexample
to the following

\medskip
\noindent {\bf Conjecture 7.1} \it The map $\Phi_h\co
T_c(S)\times\mathbb{R}_{>0}^V\rightarrow\mathbb{R}^E$ is a smooth
embedding, and the image of $\Phi_h$ is a convex polytope.\rm
\medskip

\noindent The motivation of this conjecture is as follows. Penner's simplical coordinate $\Psi$ and its deformation $\Psi_h$ are in some sense analogues to Colin de Vedi\`ere's invariant\,[CV] for  circle packings in a different setting, and the quantities $\Phi_h$ are the corresponding analogues to Rivin's invariant\,[R] for the polyhedra surfaces in this setting, see also [BS] and [L].
\medskip

\noindent 2. By Corollary 1.4, for each $h\geqslant0$, there is a
homeomorphism $$\Pi_h\co T_c(S)\times\mathbb{R}_{>0}^V\rightarrow
|A(S)-A_{\infty}(S)|\times\mathbb{R}_{>0}$$

\noindent equivariant under the mapping class group action. If
$h\neq h'$, then $\Pi_{h'}^{-1}\Pi_h$ is a self-homeomorphism of the
decorated Teichm\"{u}ller space equivariant under the mapping class
group action. These self-homeomorphisms deserve a further study. We
do not know yet if these self-homeomorphisms are smooth on the
decorated Teichm\"{u}ller space. As suggested by the referee of this article, it also seems natural to ask if these self-homeomorphisms have bounded distortion.
\\

\noindent 3. The Weil-Pertersson K\"ahler form on the Teichm\"uller space was computed in the length coordinates in \,[P2]. How to express the Weil-Petersson symplectic form on
the decorated Teichm\"{u}ller space in terms of the simplicial coordinate $\Psi$ and in terms of the $\Psi_h$ coordinate, and how to relate the $\Psi_h$ coordinate to the quantum
Teichm\"{u}ller space are interesting problems ([B], [BL], [M] and
[P3]).


\bigskip
\centerline{\bf References}
\bigskip

\noindent [BS] Bobenko, A. I.; Springborn, B. A., {\em Variational principles for circle patterns and Koebe's theorem}. Trans. Amer. Math. Soc. 356 (2004) 659–-689.
\medskip

\noindent [B] Bonahon, Francis, {\em Shearing hyperbolic surfaces,
bending pleated surfaces and Thurston's symplectic form}.  Ann. Fac.
Sci. Toulouse Math. (6)  5  (1996), 233–-297.
\medskip

\noindent [BL] Bonahon, Francis; Liu, Xiaobo, {\em Representations
of the quantum Teichm\"{u}ller space and invariants of surface
diffeomorphisms}. Geom. Topol.  11  (2007), 889–-937.
\medskip

\noindent [BE] Bowditch, Brian H.; Epstein, David B. A., {\em
Natural triangulations associated to a surface}.  Topology  27
(1988), 91–-117.
\medskip

\noindent [CV] Colin de Verdi\`ere, Yves, {\em Un principe variationnel pour les empilements de cercles}. Invent. Math. 104 (1991) 655–-669.
\medskip

\noindent [G] Guo, Ren, {\em On parameterizations of Teichm\"{u}ller
spaces of surfaces with boundary}.  J. Differential Geom.  82
(2009),  629–-640,
\medskip

\noindent [GL1] Guo, Ren; Luo, Feng, {\em Rigidity of polyhedral
surfaces. II}. Geom. Topol.  13  (2009),   1265–-1312.
\medskip

\noindent [GL2] --------,  {\em Cell decompositions of
Teichm\"{u}ller spaces of surfaces with boundary}. Preprint,
arXiv:1006.3119.
\medskip

\noindent [H] Harer, John L., {\em The virtual cohomological
dimension of the mapping class group of an orientable surface}.
Invent. Math. 84 (1986),  157–-176.
\medskip

\noindent [Ha] Hazel, Graham P., {\em Triangulating Teichm\"{u}ller
space using the Ricci flow}. Ph.D. thesis, University of California,
San Diego. 2004.
\medskip

\noindent [L] Leibon, Gregory, {\em Characterizing the Delaunay decomposition of compact hyperbolic surfaces}. Geom. Topol. 6 (2002) 361–-391.
\medskip

\noindent [L1] Luo, Feng, {\em On Teichm\"{u}ller spaces of surfaces
with boundary}.  Duke Math. J.  139  (2007),  463–-482.
\medskip

\noindent [L2] --------,  {\em Rigidity of polyhedral surfaces}.
Preprint, arXiv:math/0612714.
\medskip

\noindent [M] Mondello, Gabriele, {\em Triangulated Riemann surfaces
with boundary and the Weil-Petersson Poisson structure}. J.
Differential Geom. 81 (2009),  391–-436.
\medskip

\noindent [P1] Penner, Robert C, {\em The decorated Teichm\"{u}ller
space of punctured surfaces}.  Comm. Math. Phys. 113  (1987),
299–-339.
\medskip

\noindent[P2] --------, {\em Weil-Petersson volumes}. J. Differential Geom. 35 (1992), no. 3, 559–-608.
\medskip

\noindent [P3] --------, {\em An arithmetic problem in surface
geometry}. In: The moduli space of curves (Texel Island, 1994),
427–-466, Progr. Math., 129, Birkhauser  1995.
\medskip

\noindent [R] Rivin, Igor, {\em Euclidean structures on simplicial
and hyperbolic volume}.  Ann. of Math. (2) 139 (1994) 553–-580.
\medskip

\noindent [U] Ushiijima, Akira,  {\em A canonical cellular
decomposition of the Teichm\"{u}ller space of compact surface with
boundary}. Comm. Math. Phys. 201  (1990),   305–-326.

\bigskip
\noindent
Tian Yang\\
Department of Mathematics, Rutgers University\\
New Brunswick, NJ 08854, USA\\
(tianyang@math.rutgers.edu)

\end{document}